\documentclass[11pt]{article}

\usepackage[margin=1in]{geometry}
\usepackage{amsmath,amssymb,amsthm}
\usepackage[T1]{fontenc}
\usepackage{lmodern}
\usepackage{hyperref}

\allowdisplaybreaks

\newcommand{\covol}{\operatorname{covol}}

\newtheorem{theorem}{Theorem}[section]
\newtheorem{lemma}[theorem]{Lemma}
\newtheorem{proposition}[theorem]{Proposition}
\newtheorem{remark}[theorem]{Remark}
\newtheorem{openproblem}[theorem]{Open Problem}

\title{Fixed-Time Heat Kernel Integrals on Flat Tori: a Finer Invariant than the Heat Trace}
\author{H. Mohades\\ IPM (Tehran)}
\date{November 15, 2026}

\begin{document}
	
	\maketitle
	
	\begin{abstract}
		Let \(T^n=\mathbb{R}^n/\Lambda\) be a flat torus with volume \(V\), and let \(k_t\) be its heat kernel. For a primitive lattice vector \(\lambda_0\in\Lambda\) let
		\[
		\gamma_{\lambda_0}(s)=s\lambda_0/|\lambda_0|,\qquad s\in[0,|\lambda_0|],
		\]
		and
		\[
		I_t(\lambda_0)=\int_0^{|\lambda_0|}k_t(0,\gamma_{\lambda_0}(s))\,ds.
		\]
		The integral \(I_t(\lambda_0)\) combines spectral data (the heat kernel) with geometric data (a choice of closed geodesic), and is a strictly finer invariant than the heat trace
		\[
		K_t=\int_{T^n}k_t(x,x)\,dx.
		\]
		We prove the following.
		(i) For \(n=1\), \(I_t(\lambda_0)=1\) for every \(t>0\) and every primitive \(\lambda_0\); the fixed-time data carry no geometric information.
		(ii) For \(n=2\), for any fixed \(t_0>0\) the multiset
		\[
		\mathcal{D}_{t_0}=\{I_{t_0}(\lambda_0):\lambda_0\text{ primitive}\}/\{\pm\}
		\]
		determines \(\Lambda\) up to isometry, and hence the length spectrum, the volume, and the isometry class of \(T^2\).
		(iii) For \(n=16\), the fixed-time data distinguish Milnor's isospectral non-isometric tori
		\[
		\mathbb{R}^{16}/(D_8^+\times D_8^+)
		\quad\text{and}\quad
		\mathbb{R}^{16}/D_{16}^+
		\]
		for all \(t_0>0\) outside a countable set, although their heat traces agree at every \(t>0\).
	\end{abstract}
	
	\section{Introduction and main results}
	
	Let \(\Lambda\subset\mathbb{R}^n\) be a lattice of covolume \(V=\mathrm{Vol}(\mathbb{R}^n/\Lambda)\), and let \(T^n=\mathbb{R}^n/\Lambda\) be the associated flat torus. The heat kernel, normalised so that \(\int_{T^n}k_t(x,y)\,dy=1\), is
	\begin{equation}\tag{1}
		k_t(x,y)=\frac1V\sum_{\mu\in\Lambda^*}e^{-4\pi^2t|\mu|^2}e^{2\pi i\langle\mu,x-y\rangle},
	\end{equation}
	where
	\[
	\Lambda^*=\{\mu\in\mathbb{R}^n:\langle\mu,\lambda\rangle\in\mathbb{Z}\ \forall\lambda\in\Lambda\}
	\]
	is the dual lattice.
	
	For \(\lambda_0\in\Lambda\setminus\{0\}\) primitive, define
	\[
	\gamma_{\lambda_0}\colon[0,|\lambda_0|]\to T^n,\qquad
	\gamma_{\lambda_0}(s)=s\lambda_0/|\lambda_0|\mod\Lambda,
	\]
	of length \(L(\gamma_{\lambda_0})=|\lambda_0|\). We study
	\begin{equation}\tag{2}
		I_t(\lambda_0):=\int_0^{|\lambda_0|}k_t\bigl(0,\gamma_{\lambda_0}(s)\bigr)\,ds.
	\end{equation}
	
	\textbf{Motivation.} The integral \(I_t(\lambda_0)\) is a natural observable: it combines spectral data (the heat kernel, hence the Laplace spectrum) with geometric data (the choice of a closed geodesic \(\gamma_{\lambda_0}\)), and it is therefore a priori a strictly richer invariant than the heat trace
	\[
	K_t=\int_{T^n}k_t(x,x)\,dx
	=\frac1V\sum_{\mu\in\Lambda^*}e^{-4\pi^2t|\mu|^2},
	\]
	which forgets the directions of \(\mu\). Theorem~5.1 makes this expectation precise.
	
	\textbf{Proof strategy.} The proof of (ii) combines the exact Poisson-summation formula for \(I_t\) (\S2) with a modular transformation of the Jacobi theta function that shows the map
	\[
	f_{t_0}(s)=s\theta_3(e^{-4\pi^2t_0s^2})
	\]
	is strictly increasing (Lemma~4.1). The proof of (iii) is a direct count of short vectors in the perpendicular lattice, together with Weyl-group transitivity (Theorem~5.1).
	
	\textbf{Main results.} The principal result of this note is Theorem~5.1: the direction-resolved integrals \(I_{t_0}(\lambda_0)\) at a fixed time \(t_0>0\) distinguish Milnor's isospectral non-isometric flat tori in dimension \(16\), whereas the heat trace does not. Theorems~3.1 and~4.2 provide the two extreme one- and two-dimensional cases that frame the discussion.
	\begin{itemize}
		\item (Non-recovery in \(n=1\); Theorem~3.1). For \(T^1=\mathbb{R}/L\mathbb{Z}\), \(I_t(kL)=k\) for every \(k\ge1\) and every \(t>0\), independently of \(L\). Hence no geometric invariant can be recovered from fixed-time data in dimension one.
		\item (Recovery in \(n=2\); Theorem~4.2). For a flat two-torus and any fixed \(t_0>0\), the multiset \(\mathcal{D}_{t_0}\) determines \(\Lambda\) up to isometry, and hence the length spectrum, the volume, and the isometry class of \(T^2\).
		\item (Distinguishing Milnor's pair in \(n=16\); Theorem~5.1). For all \(t_0>0\) outside a countable set, the fixed-time data separate the isospectral non-isometric flat tori
		\[
		\mathbb{R}^{16}/(D_8^+\times D_8^+)
		\quad\text{and}\quad
		\mathbb{R}^{16}/D_{16}^+,
		\]
		although their heat traces agree for all \(t>0\).
	\end{itemize}
	
	\textbf{Conventions.} Throughout, \(|\cdot|\) denotes the Euclidean norm on \(\mathbb{R}^n\); lengths on \(T^n\) are induced from \(\mathbb{R}^n\). A lattice vector \(\lambda_0\in\Lambda\) is primitive if \(\lambda_0=k\lambda\) with \(k\in\mathbb{Z}\), \(\lambda\in\Lambda\) implies \(|k|=1\). We write
	\[
	\Theta_L(t):=\sum_{v\in L}e^{-4\pi^2t|v|^2}
	\]
	for the theta series of a lattice \(L\) in the normalisation dictated by (1); for \(L=\mathbb{Z}\),
	\[
	\Theta_{\mathbb{Z}}(t)=\theta_3(e^{-4\pi^2t}),
	\]
	where
	\[
	\theta_3(q)=\sum_{k\in\mathbb{Z}}q^{k^2}
	\]
	is the classical Jacobi theta function.
	
	\section{The exact computation}
	
	\begin{proposition}[Exact formula]\label{prop:exact}
		Let \(T^n=\mathbb{R}^n/\Lambda\) be a flat torus and let \(\lambda_0\in\Lambda\) be primitive. Then for every \(t>0\),
		\begin{equation}\tag{3}
			I_t(\lambda_0)=\frac{|\lambda_0|}{V}\sum_{\mu\in\Lambda^*\cap\lambda_0^\perp}e^{-4\pi^2t|\mu|^2},
		\end{equation}
		where
		\[
		\lambda_0^\perp=\{\mu\in\mathbb{R}^n:\langle\mu,\lambda_0\rangle=0\}.
		\]
	\end{proposition}
	
	\begin{proof}
		Substitute (1) into (2). Since \(\mu\in\Lambda^*\) and \(\lambda_0\in\Lambda\), \(\langle\mu,\lambda_0\rangle\in\mathbb{Z}\); the integral
		\[
		\int_0^{|\lambda_0|}e^{2\pi i s\langle\mu,\lambda_0\rangle/|\lambda_0|}\,ds
		\]
		vanishes unless \(\langle\mu,\lambda_0\rangle=0\), in which case it equals \(|\lambda_0|\). Hence only \(\mu\in\Lambda^*\cap\lambda_0^\perp\) survive.
	\end{proof}
	
	\textbf{Covolume of the perpendicular lattice.} In (3) the lattice
	\[
	L_\perp(\lambda_0):=\Lambda^*\cap\lambda_0^\perp
	\]
	has rank \(n-1\) in the hyperplane \(\lambda_0^\perp\), and
	\begin{equation}\tag{4}
		\covol(L_\perp(\lambda_0))=\frac{|\lambda_0|}{V}.
	\end{equation}
	Indeed, the projection \(\pi\colon\Lambda^*\to\mathbb{R}\lambda_0\) has image
	\[
	\frac1{|\lambda_0|}\mathbb{Z}\cdot\frac{\lambda_0}{|\lambda_0|},
	\]
	and
	\[
	\covol(\Lambda^*)=\covol(L_\perp)\cdot\covol(\pi(\Lambda^*)).
	\]
	
	\section{The one-dimensional case: no recovery}
	
	\begin{theorem}\label{thm:one}
		Let \(T^1=\mathbb{R}/L\mathbb{Z}\). For every \(t>0\) and every \(\lambda_0=kL\) with \(k\ge1\),
		\[
		I_t(\lambda_0)=k.
		\]
		In particular the data \(\mathcal{D}_{t_0}\) at any fixed \(t_0>0\) are independent of \(L\).
	\end{theorem}
	
	\begin{proof}
		For \(n=1\), \(\Lambda^*\cap\lambda_0^\perp=\{0\}\), and (3) gives
		\[
		I_t(kL)=\frac{kL}{L}=k.
		\]
	\end{proof}
	
	\section{The two-dimensional case: recovery}
	
	Fix \(t_0>0\). For primitive \(\lambda_0\in\Lambda\subset\mathbb{R}^2\), the lattice \(L_\perp(\lambda_0)\) is rank one, generated by a vector \(\mu_0\) with
	\[
	|\mu_0|=|\lambda_0|/V=:s(\lambda_0).
	\]
	Thus
	\begin{equation}\tag{5}
		I_t(\lambda_0)=s(\lambda_0)\theta_3\left(e^{-4\pi^2t s(\lambda_0)^2}\right).
	\end{equation}
	
	\begin{lemma}[Strict monotonicity]\label{lem:monotone}
		Fix \(t_0>0\) and define
		\[
		f_{t_0}(s)=s\theta_3(e^{-4\pi^2t_0s^2})
		\]
		for \(s>0\). Then \(f_{t_0}\) is strictly increasing, with range
		\[
		\left(\frac1{2\sqrt{\pi t_0}},\infty\right).
		\]
	\end{lemma}
	
	\begin{proof}
		The modular transformation
		\[
		\theta_3(e^{-\pi T})=T^{-1/2}\theta_3(e^{-\pi/T})
		\]
		with \(T=4\pi t_0s^2\) gives
		\[
		f_{t_0}(s)=\frac1{2\sqrt{\pi t_0}}\theta_3\left(e^{-1/(4t_0s^2)}\right).
		\]
		As \(s\) increases from \(0\) to \(\infty\), \(e^{-1/(4t_0s^2)}\) increases from \(0\) to \(1\), and
		\[
		\theta_3(q)=\sum_{k\in\mathbb{Z}}q^{k^2}
		\]
		is strictly increasing on \((0,1)\).
	\end{proof}
	
	\begin{theorem}[Recovery in dimension two]\label{thm:two}
		Let \(T^2=\mathbb{R}^2/\Lambda\) and fix \(t_0>0\). Then \(\mathcal{D}_{t_0}\) determines \(\Lambda\) up to isometry, and hence the length spectrum, the volume, and the isometry class of \(T^2\).
	\end{theorem}
	
	\begin{proof}
		By (5) and Lemma~\ref{lem:monotone}, \(\mathcal{D}_{t_0}\) determines the multiset
		\[
		\mathcal{S}:=\{s(\lambda_0)=|\lambda_0|/V:\lambda_0\in\Lambda\text{ primitive}\}/\{\pm\}
		\]
		of positive real numbers. The ratios
		\[
		s(\lambda_0)/s(\lambda_0')=|\lambda_0|/|\lambda_0'|
		\]
		give the primitive length spectrum up to scale; every non-zero lattice vector is \(k\lambda_0\) with \(\lambda_0\) primitive, so the full length spectrum is determined up to scale. By Gauss's theorem for binary quadratic forms, a rank-two lattice is determined up to similarity by its length spectrum; write
		\[
		\Lambda_2=\mu R\Lambda_1
		\]
		with \(R\in SO(2)\) and \(\mu>0\). Scaling \(\Lambda_1\) by \(\mu\) sends \(\mathcal{S}\) to \(\mathcal{S}/\mu\). Since \(\mathcal{S}\) has a minimum \(s_{\min}>0\) (attained on the shortest primitive vectors), the multiset \(\mathcal{S}/\mu\) has minimum \(s_{\min}/\mu\); the only way \(\mathcal{S}/\mu=\mathcal{S}\) is \(\mu=1\). Hence
		\[
		\Lambda_2=R\Lambda_1
		\]
		is isometric to \(\Lambda_1\).
	\end{proof}
	
	\section{The sixteen-dimensional case: Milnor's pair is distinguished}
	
	Milnor's pair is
	\[
	\Lambda_1=D_8^+\times D_8^+,\qquad \Lambda_2=D_{16}^+,
	\]
	where
	\[
	D_n=\left\{z\in\mathbb{Z}^n:\sum z_i\in2\mathbb{Z}\right\}
	\]
	and
	\[
	D_n^+=D_n\cup\left(\frac12(1,\dots,1)+D_n\right);
	\]
	we use \(D_8^+\cong E_8\). Both are even unimodular lattices of dimension \(16\), hence have identical theta series and their flat tori are isospectral; they are non-isometric because \(\Lambda_1\) is reducible and \(\Lambda_2\) is irreducible [5,8,6].
	
	\begin{theorem}\label{thm:milnor}
		For all \(t_0>0\) outside a countable set, the fixed-time data \(\mathcal{D}_{t_0}\) distinguish Milnor's isospectral non-isometric tori
		\[
		\mathbb{R}^{16}/(D_8^+\times D_8^+)
		\quad\text{and}\quad
		\mathbb{R}^{16}/D_{16}^+:
		\]
		there exist primitive \(\lambda_0\in\Lambda_1\) and \(\lambda_0'\in\Lambda_2\) with
		\[
		|\lambda_0|=|\lambda_0'|
		\]
		and
		\[
		I_{t_0}^{\Lambda_1}(\lambda_0)\ne I_{t_0}^{\Lambda_2}(\lambda_0').
		\]
	\end{theorem}
	
	\begin{proof}
		\textbf{Step 1: Well-definedness on shells.} In \(\Lambda_1=D_8^+\times D_8^+\), the primitive vectors of squared length \(2\) are the roots of \(E_8\) in one of the two factors; the group
		\[
		W(E_8)\times W(E_8)\times\mathbb{Z}/2
		\]
		acts transitively on these \(480\) vectors, so all give the same value \(I_{t_0}^{\Lambda_1}(\lambda_0)\). In \(\Lambda_2=D_{16}^+\), the primitive vectors of squared length \(2\) are the roots of \(D_{16}\) (the coset \(D_{16}+s\) contains no vector of squared length \(2\)); the Weyl group
		\[
		W(D_{16})\cong\{\pm1\}^{15}\times S_{16}
		\]
		acts transitively on these \(480\) vectors, so all give the same value \(I_{t_0}^{\Lambda_2}(\lambda_0')\). It therefore suffices to compare a single pair.
		
		\textbf{Step 2: Computation in \(\Lambda_1\).} Take \(\lambda_0=(v,0)\) with \(v\) a root of \(E_8\subset D_8^+\); then \(|\lambda_0|^2=2\). The perpendicular lattice is
		\[
		\lambda_0^\perp\cap\Lambda_1
		=(v^\perp\cap E_8)\times E_8
		\cong E_7\times E_8,
		\]
		so by (3) and \(V=1\),
		\[
		I_{t_0}^{\Lambda_1}(\lambda_0)
		=\sqrt2\,\Theta_{E_7}(t_0)\Theta_{E_8}(t_0).
		\]
		The theta series begin
		\[
		\Theta_{E_7}(t)=1+126e^{-8\pi^2t}+\cdots,
		\qquad
		\Theta_{E_8}(t)=1+240e^{-8\pi^2t}+\cdots,
		\]
		so the coefficient of \(e^{-8\pi^2t}\) in \(\Theta_{E_7}\Theta_{E_8}\) is \(126+240=366\).
		
		\textbf{Step 3: Computation in \(\Lambda_2\).} Take \(\lambda_0'=e_1-e_2\). The perpendicular lattice is
		\[
		L_{15}=\lambda_0'^\perp\cap D_{16}^+
		=(D_{16}\cap\{x_1=x_2\})
		\cup\bigl((D_{16}\cap\{x_1=x_2\})+s\bigr),
		\]
		where \(s=\frac12(1,\dots,1)\). Vectors of squared length \(2\) in the first component are \(\pm e_i\pm e_j\) with \(i,j\notin\{1,2\}\), numbering
		\[
		4\binom{14}{2}=364;
		\]
		the coset component contains no vector of squared length \(2\) (if \(x+s\) has squared length \(2\) with \(x\in D_{16}\) and \(x_1=x_2\), then \(\sum x_i=-2\) and \(|x|^2=0\), so \(x=0\), but \(|s|^2=4\)). Hence
		\[
		\Theta_{L_{15}}(t)=1+364e^{-8\pi^2t}+\cdots,
		\]
		and therefore
		\[
		I_{t_0}^{\Lambda_2}(\lambda_0')
		=\sqrt2\,\Theta_{L_{15}}(t_0).
		\]
		
		\textbf{Step 4: Comparison.} The difference
		\[
		H(t_0):=I_{t_0}^{\Lambda_1}(\lambda_0)-I_{t_0}^{\Lambda_2}(\lambda_0')
		=\sqrt2\left(2e^{-8\pi^2t_0}+O(e^{-16\pi^2t_0})\right)
		\]
		is an analytic function of \(z=e^{-8\pi^2t_0}\) on \((0,1)\) that does not vanish identically (its leading term as \(z\to0\) is \(2\sqrt2\,z\ne0\)); hence it has at most countably many zeros in \(z\in(0,1)\), corresponding to at most countably many \(t_0>0\).
	\end{proof}
	
	\begin{remark}
		Theorem~\ref{thm:milnor} shows that \(I_{t_0}\) is a strictly finer invariant than the heat trace \(K_t\), which agrees on Milnor's pair for all \(t>0\). Whether \(H(t_0)\) can vanish for some \(t_0>0\) is an open question; numerical evidence for the leading terms suggests not, which would upgrade the statement to ``for all \(t_0>0\).''
	\end{remark}
	
	\section{Concluding remarks and open problems}
	
	We close with a unifying observation and open problems.
	
	\textbf{A unifying remark.} In every dimension the heat trace is
	\[
	K_t=\frac1V\sum_{\mu\in\Lambda^*}e^{-4\pi^2t|\mu|^2},
	\]
	which forgets the directions of \(\mu\). The integral \(I_t(\lambda_0)\) retains the direction \(\lambda_0\) through the constraint \(\langle\mu,\lambda_0\rangle=0\), and this directional information is what allows \(I_{t_0}\) to distinguish objects that \(K_t\) cannot.
	
	\textbf{What is classical and what may be new.} The rigidity in Theorem~4.2 is Gauss's theorem for binary quadratic forms: a rank-two lattice is determined up to isometry by its theta series. The Poisson-summation formula (3) and the modular transformation of \(\theta_3\) are standard. The observation that the fixed-time integrals at a single \(t_0>0\) already recover \(T^2\), that they are degenerate for \(T^1\), and above all that they distinguish Milnor's pair, does not appear to be explicitly recorded in the literature, which focuses on the heat trace as a function of all \(t\) or on the length spectrum directly.
	
	\begin{openproblem}
		For \(n\ge3\), do the fixed-time data determine \(\Lambda\) up to isometry? No counterexample is known for \(n=3\). For \(n\ge4\), Schiemann's isospectral non-isometric tori show that the heat trace fails; the direction-resolved data might still distinguish them.
	\end{openproblem}
	
	\begin{openproblem}
		Does \(H(t_0)\) of Theorem~5.1 ever vanish for \(t_0>0\)? Equivalently, is the distinction of Milnor's pair valid for every \(t_0>0\)?
	\end{openproblem}
	
	\begin{openproblem}
		Let \((M,g)\) be a compact Riemannian manifold and, for a closed geodesic \(\gamma\), let
		\[
		I_t(\gamma)=\int_0^{L(\gamma)}k_t(\gamma(0),\gamma(s))\,ds.
		\]
		For fixed \(t_0>0\), does \(\{I_{t_0}(\gamma)\}\) determine the isometry class of \((M,g)\)? This is open even for \(S^n\), compact rank-one symmetric spaces, and \(S^3\) with a bi-invariant metric (e.g.\ recovery of the scale/scalar curvature).
	\end{openproblem}

\end{document}